\documentclass[12pt, leqno]{article}

\def\dim{\mathop{\rm dim}}
\def\ker{\mathop{\rm ker}}
\def\span{\mathop{\rm span}}
\def\tr{\mathop{\rm tr}}

\newtheorem{theorem}{Theorem}
\newtheorem{lemma}[theorem]{Lemma}
\newtheorem{proposition}[theorem]{Proposition}
\newtheorem{definition}[theorem]{Definition}
\newtheorem{corollary}[theorem]{Corollary}

\newcommand{\begintheorem}{\addtocounter{equation}{1}\begin{theorem}}
\newcommand{\beginlemma}{\addtocounter{equation}{1}\begin{lemma}}
\newcommand{\beginproposition}{\addtocounter{equation}{1}\begin{proposition}}
\newcommand{\begindefinition}{\addtocounter{equation}{1}\begin{definition}}
\newcommand{\begincorollary}{\addtocounter{equation}{1}\begin{corollary}}

\begin{document}

\title{Notes on groups and representations}

\author{Stephen William Semmes	\\
	Rice University		\\
	Houston, Texas}

\date{}

\maketitle


\renewcommand{\thefootnote}{}   

\footnotetext{These notes are dedicated to Bob Brooks, who told me all
sorts of cool stuff while we were visiting the Centre \'Emile Borel at
the Institut Henri Poincar\'e in the summer of 2002.}

\begin{abstract}
These informal notes are concerned with the broad themes of harmonic
analysis of groups and their representations.  We shall follow
somewhat the view of a classical analyst, with interest in various
norms in particular.  At the same time we shall try to notice some
algebraic aspects, which includes using fields other than the complex
numbers.
\end{abstract}

	Let $G$ be a \emph{group}.  Thus $G$ is a set with a
distinguished element $e$ and a binary operation, the group law, such
that $e$ is both a left and right identity element, the group
operation satisfies the associative law, and every element of $G$ has
an inverse.  If also the group operation satisfies the commutative
law, then $G$ is said to be a \emph{commutative} or \emph{abelian}
group.

	A subset $H$ of $G$ is called a \emph{subgroup} of $G$ if it
contains the identity element, if the product of any two elements of
$H$ under the group operation is also an element of $H$, and if the
inverse of each element of $H$ is also an element of $H$.  In other
words, $H$ should be a group itself using the same group operations
from $G$.

	Suppose that $G_1$, $G_2$ are groups and $\phi$ is a mapping
from $G_1$ to $G_2$.  We say that $\phi$ is a group
\emph{homomorphism} if $\phi$ maps the identity element of $G_1$ to
the identity element of $G_2$ and if $\phi$ is compatible with the
group operations on $G_1$ and $G_2$ in the sense that $\phi$ applied
to a product of elements $x$, $y$ of $G_1$ is equal to the product of
$\phi(x)$, $\phi(y)$ in $G_2$ and $\phi$ applied to the inverse of an
element $x$ of $G_1$ is equal to the inverse of $\phi(x)$ in $G_2$.

	For each subset $A$ of $G_1$, the image of $A$ under $\phi$ is
defined in the usual way as the subset of $G_2$ consisting of points
of the form $\phi(a)$, $a \in A$.  If $H$ is a subgroup of $G_1$, then
$\phi(H)$ is a subgroup of $G_2$, and in particular the image of $G_1$
under $\phi$ is a subgroup of $G_2$.

	The \emph{kernel} of $\phi$ is defined to be the subset of
$G_1$ consisting of those elements $x$ of $G_1$ with $\phi(x)$ equal
to the identity element of $G_2$.  It is easy to see that the kernel
of $\phi$ is a subgroup of $G_1$.  Also $\phi$ is \emph{injective} or
\emph{one-to-one}, meaning that $\phi$ maps $x, y \in G_1$ to the same
point in $G_2$ exactly when $x = y$, if and only if its kernel is the
trivial subgroup of the domain, consisting of the identity element
only.

	A homomorphism $\phi$ from a group $G_1$ to a group $G_2$ is
said to be an \emph{isomorphism} if $\phi$ is a one-to-one mapping
from $G_1$ onto $G_2$, which is equivalent to saying that the kernel
of $\phi$ is trivial and $\phi(G_1) = G_2$.  In this event there is an
inverse mapping $\psi$ from $G_2$ to $G_1$, characterized by the
property that $\psi(\phi(x)) = x$ for all $x \in G_1$ and
$\phi(\psi(y)) = y$ for all $y \in G_2$, and which is a group
homomorphism from $G_2$ to $G_1$.

	Let $G$ be a group and $H$ be a subgroup of $G$.  If $H$ is
the kernel of a homomorphism from $G$ to some other group, then $H$ is
a \emph{normal} subgroup of $G$, which means that $x \, h \, x^{-1}
\in H$ whenever $h \in H$ and $x \in G$.  Conversely, if $H$ is a
normal subgroup of $G$, then one can define the quotient group $G / H$
and a natural homomorphism from $G$ onto $G / H$ whose kernel is
exactly $H$.

	Now let $k$ be a \emph{field}.  This means that $k$ is a set
with two distinguished elements $0$, $1$ and two binary operations of
addition and multiplication such that $0 \ne 1$, $k$ is a commutative
group with respect to addition with $0$ as the additive identity
element, the nonzero elements of $k$ form a commutative group with
respect to multiplication with $1$ as the multiplicative identity
element, and the operations of addition and multiplication satisfy the
usual distributive laws.  This is equivalent to saying that $k$ is a
commutative ring with multiplicative identity element $1$ and that
every nonzero element of $k$ has a multiplicative inverse.  

	Recall that $k$ is said to have \emph{characteristic $0$} if
the sum of $j$ $1$'s is a nonzero element of $k$ for each positive
integer $j$.  Otherwise, there is a positive integer $j$ such that the
sum of $j$ $1$'s is equal to $0$, and the smallest such positive
integer is a prime number which is called the \emph{characteristic} of
the field $k$.

	Suppose that $V$ is a \emph{vector space} over $k$.  This
means that $V$ is a set equipped with a distinguished element $0$,
there is a binary operation on $V$ called addition so that $V$ becomes
a commutative group with additive identity element $0$, and there is
an operation of scalar multiplication which assigns to each element of
$k$ and each vector in $V$ another vector in $V$ and which enjoys
standard compatibility conditions with respect to the field operations
on $k$ and addition on $V$.  More precisely, multiplication by the
multiplicative identity element $1$ in $k$ corresponds to the identity
mapping on $V$, multiplication by any element of $k$ defines a
homomorphism on $V$ with respect to addition, etc.

	For each positive integer $n$ we get a vector space $k^n$
consisting of $n$-tuples $x = (x_1, \ldots, x_n)$ where each component
$x_l$ is an element of $k$ where the operations of addition and
multiplication by scalars are defined coordinatewise.  Namely, if $x,
y \in k^n$, then their sum $x + y$ is the element of $k^n$ whose $l$th
coordinate is given by the sum of the $l$th coordinates of $x$ and $y$
for $l = 1, \ldots, n$, and if $a \in k$ and $x \in k^n$ then the
scalar product $a \cdot x$ is the element of $k^n$ whose $l$th
coordinate is given by the product of $a$ and the $l$th coordinate of
$x$ for each $l = 1, \ldots, k$.

	A subset $L$ of a vector space $V$ over $k$ is called a
\emph{linear subspace} of $V$ if $L$ contains $0$ and if $L$ is closed
under addition and scalar multiplication.  This means that if $v, w
\in L$, then $v + w \in L$, and if $a \in k$ and $v \in L$, then the
scalar product $a \, v \in L$.  Thus $L$ is a vector space over $k$
using the operations of addition and scalar multiplication inherited
from the ones on $V$.

	If $V_1$, $V_2$ are vector spaces over the same field $k$, and
if $\phi$ is a mapping from $V_1$ to $V_2$, then we say that $\phi$ is
a \emph{linear mapping} if it is a homomorphism from $V_1$ to $V_2$ as
abelian groups, and if it preserves the operation of multiplication by
scalars in $k$.  If $\phi$ is a linear mapping from $V_1$ to $V_2$ and
if $L$ is a linear subspace of $V_1$, then the image $\phi(L)$ of $L$
under $\phi$ is a linear subspace of $V_2$.  As in the case of group
homomorphisms, the kernel of a linear mapping $\phi$ from $V_1$ to
$V_2$ is the linear subspace of $V_1$ consisting of $v \in V_1$ such
that $\phi(v)$ is the zero element of $V_2$.  The kernel of $\phi$ is
the trivial subspace of $V_1$ consisting of only the zero vector if
and only if $\phi$ is injective.  We shall write $\ker \phi$ for the
kernel of $\phi$.

	Let $V$ be a vector space over $k$, and let $v_1, \ldots, v_n$
be a finite collection of elements of $V$.  The \emph{span} of $v_1,
\ldots, v_n$, denoted $\span (v_1, \ldots, v_n)$, is the linear
subspace of $V$ consisting of all linear combinations of $v_1, \ldots,
v_n$ in $V$.  In other words the span consists of all vectors in $V$
of the form
\begin{equation}
	a_1 \, v_1 + \cdots + a_l \, v_n
\end{equation}
for some $a_1, \ldots, a_n \in k$.

	Here is another way to look at the span of $v_1, \ldots, v_l$.
Define a linear mapping $\phi$ from $k^n$ into $V$ by setting
$\phi(x)$ for $x = (x_1, \ldots, x_n) \in k^n$ to be equal to
the linear combination
\begin{equation}
	x_1 \, v_1 + \cdots + x_n \, v_n
\end{equation}
in $V$.  The span of $v_1, \ldots, v_n$ is then exactly the same as
the image of $\phi$ in $V$.

	The vectors $v_1, \ldots v_n$ in $V$ are said to be
\emph{linearly independent} if for each choice of scalars $a_1,
\ldots, a_n \in k$ we have that
\begin{equation}
	a_1 \, v_1 + \cdots + a_n \, v_n = 0
\end{equation}
if and only if the $a_j$'s are all equal to $0$.  This is equivalent
to saying that each vector in the span of $v_1, \ldots, v_n$ can be
expressed as a linear combination of $v_1, \ldots, v_n$ in a unique
way.  As another characterization, $v_1, \ldots, v_n$ are linearly
independent if and only if the linear mapping $\phi$ from $k^n$ into
$V$ defined in the previous paragraph is injective.

	Suppose that $v_1, \ldots, v_n$ and $w_1, \ldots, w_l$ are
vectors in $V$ such that $w_j \in \span (v_1, \ldots, v_n)$ for $j =
1, \ldots, l$.  A basic result in linear algebra says that if $w_1,
\ldots, w_l$ are linearly independent, then $l \le n$.  In other words,
if $l > n$, then there exist $b_1, \ldots, b_l \in k$ such that
at least one of the $b_j$'s is nonzero and
\begin{equation}
	b_1 \, w_1 + \cdots + b_l \, w_l = 0.
\end{equation}
This follows by writing the $w$'s as linear combination of $v$'s and
choosing the $b$'s so that the corresponding coefficients of the $v$'s
are all equal to $0$.  This amounts to finding a choice of $b_1,
\ldots, b_l$, not all equal to $0$, so that $n$ linear combinations of
them are equal to $0$, and this is always possible when $l > n$.

	A vector space $V$ over $k$ is said to be
\emph{finite-dimensional} if there is a finite collection of vectors
in $V$ whose span is equal to $V$.  The \emph{dimension} of $V$ is
denoted $\dim V$ and defined to be the smallest nonnegative integer
$n$ such that $V$ is the span of $n$ vectors in $V$, where the span of
$0$ vectors is defined to be simply the zero vector in $V$.

	A collection $v_1, \ldots, v_n$ of vectors in a vector space
$V$ over $k$ is said to be a \emph{basis} for $V$ if $v_1, \ldots,
v_n$ are linearly independent and if the span of $v_1, \ldots, v_n$ is
equal to $V$.  This is equivalent to saying that every vector in $V$
can be expressed in a unique way as a linear combination of $v_1,
\ldots, v_n$.  If $v_1, \ldots, v_n$ is a basis for $V$, then $V$ has
dimension equal to $n$.

	Suppose that $v_1, \ldots, v_n$ are linearly independent
vectors in a vector space $V$ of dimension $n$.  In this event
the span of $v_1, \ldots, v_n$ is equal to $V$, which is to say
that $v_1, \ldots, v_n$ is a basis for $V$.  For if $w$ is a vector
in $V$ which is not in the span of $v_1, \ldots, v_n$, then the
collection of vectors in $V$ consisting of $v_1, \ldots, v_n$
together with $w$ would also be linearly independent.

	Similarly, if $v_1, \ldots, v_n$ are vectors in a vector space
$V$ of dimension $n$ whose span is equal to $V$, then $v_1, \ldots,
v_n$ are also linearly independent and hence form a basis for $V$.
Indeed, if $v_1, \ldots, v_n$ are not linearly independent,
then one of the $v_j$'s can be expressed as a linear combination
of the others.  This would imply that $V$ is actually the span of
a proper subset of $v_1, \ldots, v_n$.

	As a basic example, let $n$ be a positive integer, and
consider the vector space $k^n$.  For $j = 1, \ldots, n$, define $e_j$
to be the vector in $k^n$ whose $j$th coordinate is equal to $1$ and
whose other coordinates are equal to $0$.  It is easy to see that
$e_1, \ldots, e_n$ form a basis for $k^n$, called the \emph{standard
basis}.

	Let $V$ be a vector space over $k$.  By a \emph{linear
functional} on $V$ we mean a linear mapping from $V$ into the scalar
field $k$, which is itself a $1$-dimensional vector space.
One can add linear functionals and multiply them by elements of $k$,
so that the space of linear functionals on $V$ is itself a vector
space over $k$.  This vector space is called the \emph{dual}
of $V$ and is denoted $V'$.

	Suppose that $v_1, \ldots, v_n$ is a basis for $V$.
For each $x \in V$ and for each integer $j$, $1 \le j \le n$,
there is a unique $\lambda_j(x) \in k$ such that
\begin{equation}
	x = \lambda_1(x) \, v_1 + \cdots + \lambda_n(x) \, v_n.
\end{equation}
In fact each $\lambda_j$ defines a linear functional on $V$.  One can
verify moreover that $\lambda_1, \ldots, \lambda_n$ forms a basis for
$V'$, called the dual basis.  In particular $V'$ also has dimension
$n$.

	If $V$, $W$ are vector spaces over $k$, then the collection of
linear mappings from $V$ to $W$ is denoted $\mathcal{L}(V, W)$.  The
dual of $V$ corresponds to the special case where $W = k$.  Just as
for the dual space, one can add linear mappings from $V$ to $W$ and
multiply them by scalars, so that $\mathcal{L}(V, W)$ becomes a vector
space over $k$ in a natural way.

	Suppose further that $v_1, \ldots, v_l$ is a basis for
$V$ and $w_1, \ldots, w_n$ is a basis for $W$.
If $j$, $m$ are integers with $1 \le j \le l$ and $1 \le m \le n$,
then define a linear mapping $\tau_{j, m}$ from $V$ to $W$
as follows.  Because $v_1, \ldots, v_l$ is a basis for $V$,
it suffices to specify the action of $\tau_{j, m}$ on this
basis, and we put
\begin{equation}
	\tau_{j, m}(v_j) = w_m, \quad
		\tau_{j, m}(v_p) = 0 \hbox{ when } p \ne j.
\end{equation}
Explicitly, $\tau_{j, m}$ sends an element
\begin{equation}
	a_1 \, v_1 + \cdots + a_l \, v_l
\end{equation}
of $V$, where $a_1, \ldots, a_l \in k$, to $a_j \, w_m \in W$.
Equivalently, if $\lambda_1, \ldots, \lambda_l$ are the linear
functionals on $V$ which are the basis for $V'$ dual to $v_1, \ldots,
v_l$, then
\begin{equation}
	\tau_{j, m}(v) = \lambda_j(v) \, w_m
\end{equation}
for each $v \in V$.  One can check that the $\tau_{j, m}$'s form a
basis for $\mathcal{L}(V, W)$.  In particular, $\mathcal{L}(V, W)$ has
dimension $l \, n$ in this situation.

	Let $V$ be a vector space over $k$, and let us write
$\mathcal{L}(V)$ for $\mathcal{L}(V, V)$, the vector space of linear
mappings from $V$ to itself.  If $T_1$, $T_2$ are linear mappings on
$V$, then we can define the product $T_1 \, T_2$ of $T_1$, $T_2$ to be
the usual composition $T_1 \circ T_2$ of $T_1$, $T_2$, which is to say
the mapping on $V$ which takes a vector $v$ to 
\begin{equation}
	T_1(T_2(v)).
\end{equation}
The product of $T_1$, $T_2$ is a linear mapping on $V$.  In this way
$\mathcal{L}(V)$ becomes an algebra over $k$, which basically means
that it is a vector space over $k$ and has this additional operation
of composition, where the composotion operation satisfies the
associative law and suitable compatibility conditions with the vector
space operations on $\mathcal{L}(V)$.  

	The identity mapping $I$ on $V$ is the mapping which takes
each vector $v \in V$ to itself, which is clearly a linear mapping.
The identity mapping on $V$ serves as an identity element in the
algebra of linear mappings on $V$ with respect to the product of
linear mappings, since the composition of $I$ with any other linear
mapping $T$ on $V$ is equal to $T$.  Thus $\mathcal{L}(V)$ is an
algebra over $k$ with an identity element.

	For each scalar $a \in k$ we get a linear mapping $a \, I$,
the product of $a$ with the identity transformation, which is the
linear mapping that sends $v$ to $a \, v$ for all $v \in V$.  The
composition of $a \, I$ with a linear transformation $T$ on $V$, in
either order, is the same as $a \, T$, which is the linear
transformation which sends $v$ to $a \, T(v)$ for each $v \in V$.

	A linear mapping $A$ on $V$ is said to be invertible if there
is another linear mapping $B$ on $V$ such that
\begin{equation}
	A \, B = B \, A = I.
\end{equation}
By standard arguments the inverse of $A$ is unique when it exists, in
which event it is denoted $A^{-1}$.  Of course the identity
transformation is invertible and is equal to its own inverse.

	If $A_1$, $A_2$ are invertible linear transformations on $V$,
then the product $A_1 \, A_2$ is invertible too, with
\begin{equation}
	(A_1 \, A_2)^{-1} = (A_2)^{-1} \, (A_1)^{-1}.
\end{equation}
If $A$ is an invertible linear transformation on $V$
and $a \in k$, $a \ne 0$, then $a \, A$ is an invertible
linear transformation on $V$ with inverse $a^{-1} \, A^{-1}$.

	A linear mapping $A$ on $V$ is invertible if it is a
one-to-one mapping of $V$ onto itself.  For in this case there is an
inverse to $A$ as a mapping from $V$ to itself, and one can easily
check that the inverse mapping is automatically linear when $A$ is.
A linear mapping on $V$ is invertible if and only if it defines
an isomorphism from $V$ onto itself.

	Suppose that $V$ is finite-dimensional, with basis
$v_1, \ldots, v_n$.  If $A$ is a one-to-one linear mapping
from $V$ into itself, then $A$ maps $v_1, \ldots, v_n$
to a set of $n$ linearly independent vectors in $V$.
Because $V$ has dimension $n$, the linearly independent
vectors $A(v_1), \ldots, A(v_n)$ span $V$, and therefore
$A$ maps $V$ onto itself.  In other words, a one-to-one
linear mapping from a finite-dimensional vector space
into itself maps the vector space onto itself and is
invertible.  If $A$, $B$ are linear transformations on $V$
such that
\begin{equation}
	B \, A = I,
\end{equation}
then $A$ is one-to-one and therefore invertible, and $B$ is the
inverse of $A$.

	Let us continue to assume that $V$ is a finite-dimensional
vector space over $k$ with basis $v_1, \ldots, v_n$.  If $A$ is a
linear mapping of $V$ onto itself, so that $A(V) = V$, then the
vectors $A(v_1), \ldots, A(v_n)$ span $V$.  Again because $V$ has
dimension $n$, it follows that $A(v_1), \ldots, A(v_n)$ are linearly
independent, and therefore form a basis of $V$.  In short, a linear
mapping from a finite-dimensional vector space onto itself is
one-to-one and hence invertible.  If $A$, $B$ are linear
transformations on $V$ such that
\begin{equation}
	A \, B = I,
\end{equation}
then $A$ maps $V$ onto itself and is therefore invertible,
and $B$ is the inverse of $A$.

	More generally, suppose that $V$, $W$ are vector spaces
over $k$ and that $T$ is a linear mapping from $V$ to $W$.
If $V$ has finite dimension, then the image of $V$ under $T$
has finite dimension, and indeed
\begin{equation}
	\dim T(V) \le \dim V.
\end{equation}
One can be more precise and say that
\begin{equation}
	\dim V = \dim (\ker T) + \dim T(V).
\end{equation}
When $V$, $W$ are finite-dimensional with the same dimension this
equation encodes the fact that the kernel of $T$ is trivial if and
only if $T$ maps $V$ onto $W$, in which case $T$ is an isomorphism
of $V$ onto $W$.

	Let $V$ be a vector space over $k$ of dimension $n$.  There is
a well-known mapping from the vector space $\mathcal{L}(V)$ of linear
mappings on $V$ into the scalar field $k$, called the
\emph{determinant}.  The determinant of a linear transformation $T$ on
$V$ is denoted $\det T$.  If we choose a basis for $V$ and express $T$
as a linear combination of the associated basis for $\mathcal{L}$, as
discussed earlier, then $\det T$ can be given explicitly as a
homogeneous polynomial of degree $n$ in the coefficients of $T$ in
this basis.

	The determinant of the identity mapping is equal to $1$, the
determinant of $a \, I$ is equal to $a^n$ for all $a \in k$, and the
determinant of a composition of linear mappings on $V$ is equal to the
product of the determinants of the individual mappings.  If $T$ is an
invertible linear transformation on $V$, then the determinant of
$T^{-1}$ times the determinant of $T$ is equal to $1$, and $\det T \ne
0$ in particular.  Cramer's rule gives an explicit formula by which
one can start with a linear transformation $A$ on $V$ and get a linear
transformation so that $A \, B$ and $B \, A$ is equal to $(\det A) I$.
When the determinant of $A$ is nonzero, it follows that $A$ is
invertible.

	Let $V$, $W$ be vector spaces over $k$, and let $V'$, $W'$ be
the corresponding dual vector spaces of linear functionals on $V$,
$W$, respectively.  If $T$ is a linear mapping from $V$ to $W$, then
we can define an associated dual mapping $T'$ from $W'$ to $V'$ as
follows.  If $\mu$ is a linear functional on $W$, then $T'(\mu)$
is the linear functional on $V$ given by
\begin{equation}
	(T'(\mu))(v) = \mu(T(v))
\end{equation}
for all $v \in V$.  Clearly $T'$ is a linear mapping from $W'$ to
$V'$.  The correspondence from $T$ to $T'$ is linear, which is to say
that it defines a linear mapping from $\mathcal{L}(V, W)$ to
$\mathcal{L}(W', V')$.

	Let us specialize now to the case where $V = W$, so that a
linear mapping $T$ on $V$ is associated to a linear mapping $T'$ on
$V'$.  Observe that the dual of the identity mapping on $V$ is the
identity mapping on $V'$.  If $T_1$, $T_2$ are linear mappings on $V$,
then
\begin{equation}
	(T_1 \, T_2)' = (T_2)' \, (T_1)',
\end{equation}
which is to say that the dual of a product of linear transformations
is equal to the product of the corresponding dual linear
transformations in the opposite order.

	In general, if $T$ is an isomorphism from the vector
space $V$ onto the vector space $W$, then the dual $T'$ is an
isomorphism from the dual space $W'$ onto the dual space $V'$.
If $V = W$ and $T$ is an invertible linear transformation
from $V$ to itself, then the remarks in the preceding paragraph
show that the dual of $T^{-1}$ is the same as the inverse
of the dual transformation $T'$.

	Let $V$ be a vector space over $k$, and let us write $GL(V)$
for the group of invertible linear transformations on $V$, using
composition of linear mappings as the group operation.  This is called
the general linear group associated to $V$.  The linear
transformations of the form $a \, I$, $a \in k$, $a \ne 0$, form a
commutative subgroup of $GL(V)$, and when $\dim V = 1$ this is all of
$GL(V)$.  The determinant defines a homomorphism from $GL(V)$ to the
commutative multiplicative group of nonzero elements of $k$.

	Let $V$ be a finite-dimensional vector space over $k$ with
basis $v_1, \ldots, v_n$, and let $T$ be a linear transformation from
$T$ to itself.  If $j$, $l$ are integers with $1 \le j, l \le n$, then
we have the linear mapping $\tau_{j, l}$ on $V$ as before, with
$\tau_{j, l}(v_j) = v_l$ and $\tau_{j, l}(v_p) = 0$ when $p \ne j$.
These linear transformations $\tau_{j, l}$, $1 \le j, l \le n$, form a
basis for the vector space of linear transformations on $V$, as
discussed previously.  Thus $T$ can be expressed in a unique manner as
a linear combination of the $\tau_{j, l}$'s.

	The \emph{trace} of $T$, denoted $\tr T$, is defined to be the
sum of the $\tau_{j, j}$ coefficients of $T$, $1 \le j \le n$, in the
expansion of $T$ as a linear combination of $\tau_{j, l}$'s.  It
follows that $\tr T$ is a linear function of $T$, i.e., a linear
functional on the vector space of linear mappings on $V$.  Notice that
the trace of the identity mapping is equal to the sum of $n$ $1$'s in
$k$, where $n$ is the dimension of $V$.

	A basic property of the trace states that if $A$, $B$
are linear mappings on $V$, then
\begin{equation}
	\tr A \, B = \tr B \, A.
\end{equation}
This can be computed directly from the definition.  In particular, if
$A$, $T$ are linear transformations on $V$ and $A$ is invertible, then
\begin{equation}
	\tr A^{-1} \, T \, A = \tr A.
\end{equation}
As a consequence one can show that the trace of a linear
transformation does not depend on the choice of basis for $V$.  For if
one had a second basis for $V$, then one could pass from the first
choice of basis to the second one using an invertible linear
transformation, and the preceding identity implies that the
definitions of the trace associated to the two bases coincide.

	Let us now consider some aspects of functional calculus.  Let
$k$ be a field, and let us write $\mathcal{P}(k)$ for the polynomial
algebra with coefficients in $k$.  Thus an element $p(t)$ of
$\mathcal{P}(k)$ is given by a formal sum of the form
\begin{equation}
	p(t) = c_l \, t^l + c_{l-1} \, t^{l - 1} + \cdots + c_0,
\end{equation}
where $l$ is a nonnegative integer and $c_0, \ldots, c_l \in k$.  In
this case we way that $p(t)$ has degree less than or equal to $l$, or
equal to $l$ if $c_l \ne = 0$.

	More precisely, an element $p(t)$ of $\mathcal{P}(k)$
determines a function on $k$ in a natural way, but we think of $p(t)$
as being more specific than that.  The coefficients $c_0, \ldots, c_l$
are part of the data, although of course one can add terms with
coefficients equal to $0$ without changing the element of
$\mathcal{P}(k)$.  We shall say more about this in a moment.

	Elements of $\mathcal{P}(k)$ can be added in the usual manner,
term by term.  One can also multiply an element of $\mathcal{P}(k)$ by
an element of $k$, so that $\mathcal{P}(k)$ is an infinite-dimensional
vector space over $k$.  Moreover, one can multiply two elements of
$\mathcal{P}(k)$, so that $\mathcal{P}(k)$ is a commutative algebra
over $k$.  The polynomial algebra $\mathcal{P}(k)$ contains a copy of
$k$ as constant polynomials, and the constant $1$ is the
multiplicative identity element of $\mathcal{P}(k)$.

	If $p(t)$ is an element of $\mathcal{P}(k)$, then the function
on $k$ associated to $p$ vanishes at $0$ if and only if the constant
term in $p(t)$ is equal to $0$, and this is equivalent to saying that
$p(t)$ can be expressed as $t \, q(t)$ for some $q(t) \in
\mathcal{P}(k)$.  More generally, the function on $k$ associated to
$p(t) \in \mathcal{P}(k)$ vanishes at some $a \in k$ if and only if
$p(t)$ can be expressed as $(t - a) \, q(t)$ for some $q(t) \in
\mathcal{P}(k)$.  The function on $k$ associated to $p(t)$ vanishes at
the distinct points $a_1, \ldots, a_n \in k$ if and only if $p(t)$ can
be expressed as
\begin{equation}
	(t - a_1) \cdots (t - a_n) \, q(t)
\end{equation}
for some $q(t) \in \mathcal{P}(k)$.  In particular, if $p(t)$ is not
the zero polynomial, which is to say that $p(t)$ has at least
one nonzero coefficient, and if $p(t)$ has degree less than or equal
to $l$, then the function on $k$ associated to $p(t)$ can vanish
on a subset of $k$ with at most $l$ elements.

	If $k$ is infinite, then it follows that the function on $k$
associated to $p(t)$ is equal to $0$ at every point in $k$ if and only
if all the coefficients of $p(t)$ are equal to $0$.  If $k$ has
characteristic $0$, then $k$ contains a copy of the rational numbers,
and therefore $k$ is infinite.  If $k$ has positive characteristic
then $k$ may or may not be finite.

	If $k$ is finite, then there certainly are elements $p(t)$ of
$\mathcal{P}(k)$ for which the associated function on $k$ vanishes at
every element of $k$, even though $p(t)$ has nonzero coefficients and
is therefore not equal to $0$ in $\mathcal{P}(k)$.  An element $p(t)$
of $\mathcal{P}(k)$ can just as well be viewed as an element of the
polynomial algebra with coefficients in any field which contains $k$,
and hence defines a function on any field which contains $k$.  To say
that $p(t)$ is equal to $0$ as an element of $\mathcal{P}(k)$, which
means that all of its coefficients are equal to $0$, is equivalent to
saying that the function associated to $p(t)$ on any extension of $k$
to a larger field is equal to $0$ at every point in the field.  In
other words, if $p(t)$ has a nonzero coefficient, then the function
associated to $p(t)$ on some extension of $k$ is different from $0$ at
some point in the larger field.

	We can extend this further by letting elements of
$\mathcal{P}(k)$ operate on linear transformations.
Let $V$ be a vector space over $k$, and let $A$ be a linear
transformation on $V$.  If
\begin{equation}
	p(t) = c_l \, t^l + \cdots c_0
\end{equation}
is an element of $\mathcal{P}(k)$, so that $c_0, \ldots, c_l \in k$,
then we can define $p(A)$ to be the linear transformation on $V$ given
by
\begin{equation}
	p(A) = c_l \, A^l + \cdots c_0 \, I.
\end{equation}
Here $A^j$ is the $j$th power of $A$ for each positive integer $j$,
which means that $A^j$ is the product of $j$ copies of $A$.

	Let $A$ be a linear transformation on $V$, and let
$p_1(t)$, $p_2(t)$ be elements of $\mathcal{P}(k)$.
The sum $p_1 + p_2$ and product $p_1 \, p_2$ are also elements
of $\mathcal{P}(k)$, and it is easy to see that
\begin{equation}
	(p_1 + p_2)(A) = p_1(A) + p_2(A), \quad
		(p_1 \, p_2)(A) = p_1(A) \, p_2(A).
\end{equation}
In other words the action of $\mathcal{P}(k)$ on $\mathcal{L}(V)$ is
compatible with the operations of addition and multiplication in the
obvious manner.  Notice that if $A$ is of the form $a \, I$, where $a
\in k$, then $p(A)$ is equal to $p(a) \, I$, where $p(a)$ is the value
of the function on $k$ associated to $p(t)$ at $a$.

	Let $G$ be a finite group with $n$ elements, and let $k$ be a
field.  By a \emph{representation} of $G$ over $k$ we mean a pair
$(\rho, V)$, where $V$ is a vector space over $k$ of positive and
finite dimension and $\rho$ is a homomorphism from $G$ into the group
$GL(V)$ of invertible linear transformations on $V$.  Explicitly, for
each $x \in G$ we shall write $\rho_x$ for the corresponding linear
transformation on $V$, so that if $v \in V$ then $\rho_x(v)$ is its
image under $\rho_x$.  To say that $\rho$ is a homomorphism from $G$
into $GL(V)$ means that $\rho_e = I$, where $e$ is the identity
element of $G$,
\begin{equation}
	\rho_{x \, y} = \rho_x \, \rho_y
\end{equation}
for all $x, y \in G$, and that
\begin{equation}
	\rho_{x^{-1}} = (\rho_x)^{-1}
\end{equation}
for all $x \in G$.

	The \emph{degree} of a representation of $G$ is defined to be
the dimension of the vector space on which the representation acts.
Two representations $(\rho, V)$ and $(\sigma, W)$ of $G$ are said to
be \emph{isomorphic} if there is a one-to-one linear mapping $\phi$
from $V$ onto $W$ such that
\begin{equation}
	\phi \circ \rho_x = \sigma_x \circ \phi
\end{equation}
for all $x \in G$.  In other words, $\phi$ should be a linear
isomorphism from $V$ onto $W$ which intertwines the representations
$\rho$, $\sigma$.  Note that isomorphic representations have the same
degree.

	Let $(\rho, V)$ be a representation of $G$.  Define a function
$\lambda$ on $G$ associated to this representation by
\begin{equation}
	\lambda(x) = \tr \rho_x,
\end{equation}
i.e., $\lambda(x)$ is the trace of the linear transformation $\rho_x$
on $G$.  This is the \emph{character} associated to $(\rho, V)$, and
it is a function on $G$ with values in $k$.  At the identity element
$e$ of $G$ the value of $\lambda$ is equal to the sum of $1$'s in $k$
where the number of $1$'s is the degree of the representation.

	If $G$ is any group, then two elements $x$, $y$ of $G$
are said to be \emph{conjugate} if there is a $w \in G$ such that
\begin{equation}
	y = w \, x \, w^{-1}.
\end{equation}
Clearly an element of $G$ is conjugate to itself, and conjugacy is
also symmetric in the two group elements.  It is transitive as well,
which means that if $x, y, z \in G$, $x$ is conjugate to $y$, and $y$
is conjugate to $z$, then $x$ is conjugate to $z$.  In short,
conjugacy defines an \emph{equivalence relation} on $G$.  Thus the
group $G$ can be partitioned into equivalence classes, called
conjugacy classes, where two elements of $G$ lie in the same
equivalence class exactly when they are conjugate.

	A function on $G$ is called a class function if it is constant
on conjugacy classes.  If $\lambda$ is the character associated to a
representation $(\rho, V)$ of $G$, then $\lambda$ is a class function,
because
\begin{equation}
	\lambda(w \, x \, w^{-1}) = \tr (\rho_{w \, x \, w^{-1}})
		= \tr (\rho_w \, \rho_x \, (\rho_w)^{-1})
		= \tr \rho_x = \lambda(x)
\end{equation}
for all $x, w \in G$.  For similar reasons notice that the characters
associated to two isomorphic representations are equal to each other.
More precisely, if $V$ and $W$ are finite-dimensional vector spaces
over the same field $k$, if $T$ is a linear mapping on $V$, and if
$\phi$ is a linear isomorphism of $V$ onto $W$, then the trace of $T$
as a linear mapping on $V$ is equal to the trace of $\phi \circ T
\circ \phi^{-1}$ as a linear mapping on $W$.  A closely related
statement is that if $A$ is a linear mapping from $V$ to $W$ and $B$
is a linear mapping from $W$ to $V$, then the trace of $B \circ A$ as
a linear transformation on $V$ is equal to the trace of $A \circ B$ as
a linear mapping on $W$.

	Notice that a homomorphism from a group $G$ into an abelian
group is automatically a class function.  Also, a representation of a
finite group $G$ of degree $1$ is basically a homomorphism of $G$ into
an abelian group, namely, a homomorphism into the multiplicative group
of nonzero elements of $k$.  The character of a representation of
degree $1$ exactly gives this homomorphism.  For if $V$ is a vector
space over $k$ of dimension $1$ and $T$ is a linear transformation on
$V$, then there is an $a \in k$ such that $T(v) = a \, v$ for all $v
\in V$, and the trace of $T$ is exactly equal to $a$.

	If $V$ is a vector space over $k$ of positive finite dimension
$\ell$, then $V$ is isomorphic as a vector space over $k$ to $k^\ell$.
Thus every representation of a finite group $G$ over the field
$k$ is isomorphic to a representation on $k^\ell$ for some positive
integer $\ell$.  Using the standard basis for $k^\ell$, the
general linear group over $k^\ell$ can be identified with the group
of $\ell \times \ell$ invertible matrices with entries in $k$.
For that matter one could start with an $\ell$-dimensional
vector space $V$, choose a basis for $V$, and use that to
identify linear transformations on $V$ with $\ell \times \ell$
matrices with entries in $k$.

	At any rate, in general a representation of a finite group $G$
over a field $k$ is basically the same thing as a homomorphism from
$G$ into the group of $\ell \times \ell$ invertible matrices with
entries in $k$.  To get the character of the representation one
takes the traces of the corresponding matrices, which is to say the
sum of the diagonal entries.  When $\ell = 1$ the matrix and the
trace are basically the same thing, an element of $k$.

	Part of the business with representations is that one can
mess with the field $k$, and this is indeed a fascinating matter.
If $k$ is a subfield of a larger field, then a representation of
a finite group $G$ over $k$ leads to a representation of $G$
over the larger field in a natural way.  This is especially clear
in terms of matrices, because matrices with entries in $k$ can also
be viewed as matrices with entries in a larger field.
Even though the new representation is viewed as acting on vector
spaces over the larger field, the character obtained in this way
will be the same as the character of the original representation
over $k$, and in particular it still takes values in $k$.

	One can look at this in the other direction and start with a
representation of $G$ over $k$, and ask if it perhaps comes from a
representation over a subfield of $k$.  A necessary condition for this
to happen is that the character should take values in the subfield.
It may be that the representation is described initially in terms of
matrices with entries in $k$ which are not contained entirely in the
subfield, but that an isomorphic realization uses only matrices with
entries in the subfield.

	There is a way to restrict to a proper subfield that works
automatically.  Namely, if $k$ is a field and $V$ is a vector space
over $k$, then $V$ is also a vector space over any subfield of $k$.
A linear mapping on $V$ with respect to $k$ is also linear with
respect to any subfield of $k$.  In this way a representation
over $k$ can be converted into a representation over a subfield,
with a suitable increase in the degree of the representation.

	Another basic scenario is that one has a representation of a
finite group $G$ over the rational numbers, say, and that the
representation can be described in terms of matrices with integer
entries.  In this event one can try to reduce modulo $p$ to get a
representation of $G$ over the the field ${\bf Z}_p$ of integers
modulo $p$, where $p$ is a prime number.  Conversely one might start
with a representation over ${\bf Z}_p$ and ask whether it arises from
reduction modulo $p$ of a representation over the rational numbers.
Of course there are a lot of variations of these themes.

	Let us now consider some basic examples of representations.
Fix a finite group $G$ with $n$ elements and a field $k$.  One
automatically has the unit representation on the one-dimensional
vector space $k$, in which every element of $G$ is sent to the
identity transformation on $k$.  The character associated to this
representation is equal to $1$ at each point in $G$.

	If $A$ is a finite nonempty set, then let us write
$\mathcal{F}(A, k)$ for the vector space of $k$-valued functions on
$X$.  The dimension of this vector space is equal to the number of
elements of $A$.  Suppose that we have an action of $G$ on $A$, which
means a homomorphism from $G$ into the group of permutations on $A$.
In other words, suppose that for each $x \in G$ we have a
one-to-one mapping $\pi_x$ from $A$ onto itself, which is to say a
permutation on $A$, such that $\pi_e$ is the identity mapping on $A$,
\begin{equation}
	\pi_{x \, y} = \pi_x \circ \pi_y
\end{equation}
for all $x, y \in G$, and
\begin{equation}
	\pi_{x^{-1}} = (\pi_x)^{-1}
\end{equation}
for all $x \in G$.  This leads to a representation of $G$ on
$\mathcal{F}(A, k)$, by composing functions on $A$ with these
permutations in an appropriate manner.  Specifically, for each $x \in
G$, we use the linear transformation on $\mathcal{F}(A, k)$ which
takes a function $f$ to $f \circ (\pi_x)^{-1}$, where the inverse is
employed so that the composition laws come out in the right order.

	Here is another way to look at this representation.  For each
$a \in A$, let $\delta_a$ be the function on $\mathcal{F}(A, k)$ which
is equal to $1$ at $a$ and to $0$ at all other elements of $A$.  The
collection of functions $\delta_a$, $a \in A$, form a basis for
$\mathcal{F}(A, k)$.  If $x \in G$, then the representation of $G$ on
$\mathcal{F}(A, k)$ just described is characterized by the fact that
the linear transformation on $\mathcal{F}(A, k)$ associated to $x$
sends $\delta_a$ to $\delta_b$ with $b = \pi_x(a)$.  The character
associated to this representation at a point $x$ in $G$ is equal to a
sum of $1$'s, where the number of $1$'s is the number of fixed points
of $\pi_x$ on $A$.

	For instance, one can take $A = G$ and define $\pi_x$ to be
the permutation on $G$ given by left multiplication by $x$.  This
leads to the left regular representation of $G$.  Instead one can take
$\pi_x$ to be right multiplication by $x^{-1}$, and this leads to the
right regular representation of $G$.  These representations are
isomorphic to each other, as one can see by using the mapping $y
\mapsto y^{-1}$ on $G$ to switch from one action to another.  The
character of the regular representation is given by the function on
$G$ equal to the sum of $n$ $1$'s at the identity element of $G$ and
equal to $0$ at all other elements of $G$.

	More generally, suppose that $H$ is a subgroup of $G$.  One
can then define the space of cosets $G / H$ in the usual manner, and
this space comes equipped with an action by $G$ which leads to a
representation of $G$.  The subgroup $H$ is not required to be a
normal subgroup; that would be needed in order for $G / H$ to be a
group, but one can define the quotient as a set with a $G$ action for
any subgroup $H$.

	Let $k$ be a field, and let $V$, $W$ be vector spaces over
$k$.  One can define the direct sum of $V$ and $W$ in such a way that
the direct sum is a vector space over $k$ containing copies of $V$,
$W$, and in which every element of the direct sum can be expressed as
a sum of elements in the copies of $V$ and $W$ in a unique manner.  If
$V$ and $W$ are finite-dimensional, then the direct sum is also finite
dimensional, with dimension equal to the sum of the dimensions of $V$
and $W$.  If $A$, $B$ are linear transformations on $V$, $W$, then
there is a linear transformation on the direct sum which maps the
copies of $V$ and $W$ to themselves and whose restrictions to the
copies of $V$ and $W$ in the direct sum are equal to $A$, $B$.  If
$V$, $W$ have finite dimension, then the trace of the combined linear
transformation on the direct sum is equal to the sum of the traces of
$A$ and $B$ on $V$ and $W$, respectively.

	Suppose that $G$ is a finite group, $k$ is a field, and
$(\rho, V)$ and $(\sigma, W)$ are representations of $G$ over $k$.
There is a natural way to take the direct sum of these two
representations, acting on the direct sum of $V$ and $W$.  Namely, for
each $x \in G$, we have the linear transformations $\rho_x$ on $V$ and
$\sigma_x$ on $W$, and these can be combined to give a linear
transformation on the direct sum as in the preceding paragraph.
The character of the direct sum representation is equal to the sum of
the characters associated to $\rho$ and $\sigma$.

	If $V$ and $W$ are vector spaces over a field $k$, then there
is a standard construction of a tensor product vector space over $k$.
If $V$ and $W$ have finite dimension, then so does the tensor product,
and the dimension of the tensor product is equal to the product of the
dimensions of $V$ and $W$.  If $G$ is a finite group, $k$ is a field,
and $(\rho, V)$ and $(\sigma, W)$ are representations of $G$, then we
get a tensor product representation acting on the tensor product of
$V$ and $W$.  The character of the tensor product representation is
equal to the product of the characters associated to $\rho$, $\sigma$.

	Let $G$ be a finite group, let $k$ be a field, and let $(\rho,
V)$ be a representation of $G$.  As discussed earlier we can define
the dual vector space $V'$ consisting of the linear functionals on
$V$, which has the same dimension as $V$, and each linear
transformation $T$ on $V$ leads to a dual linear transformation $T'$
on $V'$.  The representation dual to $(\rho, V)$ acts on $V'$ by
sending $x \in G$ to the dual of $(\rho_x)^{-1}$.  If $\lambda(x)$ is
the character associated to $\rho$, then the character associated to
the dual representation is given by $\lambda(x^{-1})$.  This uses the
fact that if $T$ is a linear transformation on a finite-dimensional
vector space $V$, then the trace of the dual linear transformation
$T'$ on $V'$ is equal to the trace of $T$ on $V$.

	Let $G$ be a finite group with $n$ elements, let $k$ be a
field, and let $\mathcal{F}(G, k)$ denote the vector space of
$k$-valued functions on $G$ as before.  For each $x \in G$ we again
write $\delta_x$ for the function on $G$ which is equal to $1$ at $x$
and to $0$ at all other elements of $G$.  This is a basis for
$\mathcal{F}(G, k)$, which has dimension $n$ as a vector space
over $k$.

	Suppose that $f_1$, $f_2$ are $k$-valued functions on $G$.
The \emph{convolution} of $f_1$, $f_2$ is the $k$-valued function on
$G$ denoted $f_1 * f_2$ and given by
\begin{equation}
	(f_1 * f_2)(z) 
		= \sum_{x, y \in G \atop z = x \, y} f_1(x) \, f_2(y).
\end{equation}
This operation of convolution is associative and satisfies the
distributive laws with respect to addition and scalar multiplication,
which is to say that it is linear in $f_1$, $f_2$.  Thus the vector
space $\mathcal{F}(G, k)$ becomes an algebra.

	For each $x, y \in G$ we have that the convolution of
$\delta_x$ and $\delta_y$ is equal to $\delta_z$, with $z = x \, y$.
In other words, on the $\delta_x$'s, the convolution reduces exactly
to the group operation on $G$.  Convolution of arbitrary functions on
$G$ is determined by this and linearity, since the $\delta_x$'s form a
basis for $\mathcal{F}(G, k)$.  If $e$ is the identity element of $G$,
then convolution of any function $f$ on $G$ with $\delta_e$ gives $f$
back again, which is to say that $\delta_e$ is the identity element in
$\mathcal{F}(G, k)$ for convolution.  Convolution with other
$\delta_x$'s is given by translation of the function.

	Let us write $\mathcal{F}_c(G, k)$ for the subspace of
functions on $G$ which are class functions, i.e., which are constant
on the conjugacy classes of $G$.  Thus the dimension of
$\mathcal{F}_c(G, k)$ is equal to the number of conjugacy classes in
$G$.  One can also characterize $\mathcal{F}_c(G, k)$ as the center of
the convolution algebra $\mathcal{F}(G, k)$.  Namely, a function on
$G$ is a class function if and only if it commutes with all other
functions on $G$ with respect to convolution.  This is equivalent to
saying that it commutes with all of the $\delta_x$'s, $x \in G$.

	Suppose that $(\rho, V)$ is a representation of $G$ over $k$.
If $f \in \mathcal{F}(G, k)$, then we can associate a linear
transformation $T_f$ on $V$ to $f$ using the representation, namely,
\begin{equation}
	T_f = \sum_{x \in G} f(x) \, \rho_x.
\end{equation}
The correspondence $f \mapsto T_f$ is clearly linear, and when $f =
\delta_x$ for some $x \in G$ we have that $T_f$ is equal to $\rho_x$.
In fact the correspondence $f \mapsto T_f$ is an algebra homomorphism,
which is to say that the convolution of two functions is sent to the
composition of the associated linear transformations on $V$.

	Let $G$ be a finite group with $n$ elements, let $k$
be a field, and let $(\rho, V)$ be a representation of $G$.
A linear subspace $L$ of $V$ is said to be invariant under the
representation if
\begin{equation}
	\rho_x(L) \subseteq L
\end{equation}
for all $x \in G$, which is equivalent to
\begin{equation}
	\rho_x(L) = L
\end{equation}
for all $x \in G$.  We say that $(\rho, V)$ is \emph{irreducible} if
the only linear subspaces $L$ of $V$ which are invariant under the
representation are the trivial subspace consisting of only the zero
vector and $V$ itself.  

	One-dimensional representations are automatically irreducible.
In general a representation $(\rho, V)$ of $G$ is irreducible if
and only if for each $v \in V$ with $v \ne 0$ we have that
\begin{equation}
	V = \span \{\rho_x(v) : x \in G\}.
\end{equation}
Indeed, the span of the vectors $\rho_x(v)$, $x \in G$, is automatically
invariant under the representation, and so must be all of $V$ is the
representation is irreducible.  Conversely, if $L$ is a linear subspace
of $V$ which is invariant under the representation and which contains
a nonzero vector $v$, then $L$ contains the vectors $\rho_x(v)$ for
$x \in G$, and therefore contains their span, which is all of $V$
by assumption.

	Suppose that $(\rho, V)$ and $(\sigma, W)$ are representations
of $G$.  Let $\phi$ be a linear mapping from $V$ to $W$ which
intertwines the representations, which is to say that
\begin{equation}
	\sigma_x \circ \phi = \phi \circ \rho_x
\end{equation}
for all $x \in G$.  The kernel and image of $\phi$ are linear
subspaces of $V$ and $W$ which are invariant under the representations
$\rho$, $\sigma$, as one can easily check.

	If $(\rho, V)$ is irreducible, then $\phi$ must either be the
zero mapping or injective.  If $(\sigma, W)$ is irreducible, then
$\phi$ is either the zero mapping or it maps $V$ onto $W$.  If both
$(\rho, V)$ and $(\sigma, W)$ are irreducible, then $\phi$ is either
the zero mapping or an isomorphism.  In particular, either
$\phi$ is the zero mapping, or $(\rho, V)$ and $(\sigma, W)$ are
isomorphic representations, and thus have the same character.
These statements constitute one-half of ``Schur's lemma''.

	Let $V$ be a finite-dimensional vector space over a field $k$,
let $T$ be a linear tranformation on $V$, and let $L$ be a nonzero
proper linear subspace of $V$ which is invariant under $T$, so that
$T(L) \subseteq L$.  Of course one can restrict $T$ to $L$ to get a
linear transformation there.  One can also form the quotient $V / L$,
and $T$ induces a linear transformation on the quotient.  The trace of
$T$ as a linear transformation on $V$ is equal to the sum of the trace
of the restriction of $L$ to $T$ and the trace of the linear mapping
on $V / L$ induced by $T$.  It may or may not be that there is a
linear subspace of $V$ complementary to $L$ which is invariant under $L$,
which is to say that $V$ would be isomorphic to a direct sum of two
vector spaces in such a way that $T$ would correspond to a sum of
two linear operators on each of the two pieces separately.

	Suppose that $G$ is a finite group with $n$ elements, $k$ is a
field, and $(\rho, V)$ is a representation of $G$ over $k$.  Suppose
further that $L$ is a proper nonzero linear subspace of $V$ which is
invariant under the representation.  Thus we can restrict the
representation to $L$ to get a new representation of $G$.  We can also
form the quotient space $V / L$, and the representation on $V$ induces
one on $V / L$ because $L$ is invariant.  The character associated to
the original representation on $V$ is equal to the sum of the
characters associated to the restriction of the representation to $L$
and to the representation induced on the quotient $V / L$.

	It may or may not be that the representation $(\rho, V)$ is
actually isomorphic to the direct sum of these two representations of
smaller degree.  This amounts to the question of whether there is a
linear subspace of $V$ which is complementary to $L$ and invariant
under the representation.  It turns out that this does always happen
if either $k$ has characteristic $0$, or if $k$ has positive
characteristic $p$ and the order $n$ of $G$ is not an integer multiple
of $p$.

	Recall that a projection of $V$ onto $L$ is a linear mapping
on $V$ which sends every vector in $V$ into $L$, and which sends
every vector in $L$ to itself.  The kernel of the projection is
a linear subspace of $V$ which is complementary to $L$, and conversely
one can start with a linear subspace of $V$ complementary to $L$ and
get a projection of $V$ onto $L$ with that subspace as its kernel.
The question of having an invariant complement to $L$ is equivalent
to having a projection of $V$ onto $L$ which commutes with the
representation.

	Suppose that $P$ is any projection of $V$ onto $L$.  For each
$x \in G$, $\rho_x \circ P \circ (\rho_x)^{-1}$ is another projection
on $V$, and the image of this projection is also equal to $L$ because
$L$ is invariant under the representation.  The conditions above on
$k$ are equivalent to saying that a sum of $n$ $1$'s in $k$ is not
equal to $0$ in $k$.  This permits one to average over $x \in G$ to
get a projection of $V$ onto $L$ which commutes with the
representation by construction.

	From now on in these notes let us assume that
\begin{equation}
	\hbox{$k$ has characteristic equal to $0$}.
\end{equation}
This implies that every representation of a finite group $G$ over $k$
is isomorphic to a direct sum of irreducible representations.
Otherwise one would get a kind of composition series of irreducible
representations.  In any event every character of a representation of
$G$ over $k$ can be expressed as a sum of characters of irreducible
representations.

	Let us assume further that
\begin{equation}
	\hbox{$k$ is algebraically closed},
\end{equation}
which means that every nonconstant polynomial on $k$ has a root, and
hence can be factored.

	Let $V$ be a finite-dimensional vector space over $k$, and let
$T$ be a linear transformation on $V$.  An interesting polynomial
associated to $T$ is the characteristic polynomial $p(\alpha) = \det
(T - \alpha \, I)$.  Because $k$ is algebraically closed, there is a
$\alpha \in k$ such that $p(\alpha) = 0$.  This is equivalent to
saying that there is a $\alpha \in k$ such that $T - \alpha \, I$ is
not invertible.

	Let $V$ be a finite-dimensional vector space over $k$, let $T$
be a linear operator on $V$, and let $\alpha$ be an element of $k$.
A linear operator on $V$ is invertible if and only if it has trivial
kernel, and thus $T - \alpha \, I$ is not invertible if and only if
$T - \alpha \, I$ has a nontrivial kernel.  The kernel of $T -
\alpha \, I$ will be denoted $E(T, \alpha)$ and a vector $v \in V$
lies in $E(T, \alpha)$ if and only if $T(v) = \alpha \, v$.  As in
the preceding paragraph, for each linear transformation $T$ on $V$
there is a $\alpha \in k$ such that $E(T, \alpha)$ is nontrivial.
When $E(T, \alpha)$ is nontrivial, which is to say that it contains
nonzero vectors, then we say that $\alpha$ is an eigenvalue of $T$,
and $E(T, \alpha)$ is the corresponding eigenspace of eigenvectors of
$T$ with eigenvalue $\alpha$.

	Let $V$ be a finite-dimensional vector space over $k$, let
$A$, $T$ be linear transformations on $V$, and let $\alpha \in k$ be
an eigenvalue for $T$.  Suppose that $A$ and $T$ commute, which is to
say that $A \, T = T \, A$.  If $v \in V$ is an eigenvector for $T$
with eigenvalue $\alpha$, so that $T(v) = \alpha \, v$, then $A(v)$
is too.  In other words, the eigenspace $E(T, \alpha)$ is invariant
under $A$.

	Let $G$ be a finite group with $n$ elements, and let $(\rho,
V)$ be an irreducible representation of $G$ over $k$.  Suppose that
$T$ is a linear transformation on $V$ which commutes with the
representation, which is to say that $T \circ \rho_x = \rho_x \circ T$
for all $x \in G$.  If $\alpha \in k$ is an eigenvalue of $T$, then
the corresponding eigenspace $E(T, \alpha)$ is a nonzero linear
subspace of $V$ which is invariant under the representation.
Irreducibility implies that $E(T, \alpha) = V$, which is to say that
$T = \alpha \, I$.  This is the second part of Schur's lemma.

	Assume also that $f$ is a $k$-valued class function on $G$.
As before let $T_f$ be the linear transformation on $V$ given by
$\sum_{x \in G} f(x) \, \rho_x$.  The assumption that $f$ is a class
function implies that $T_f$ commutes with the representation $\rho$.
Thus $T_f$ is equal to a scalar multiple of the identity
transformation on $V$, by Schur's lemma.

	If $G$ happens to be an abelian group, and $(\rho, V)$ is a
representation over $G$, then $\rho_y$ commutes with the
representation for all $y \in G$.  If the representation is
irreducible, then it follows that $\rho_y$ is a scalar multiple of the
identity for all $y \in G$.  In fact the representation has degree
equal to $1$ in this case.  In other words, the representation is
given by a homomorphism from $G$ into the multiplicative group of
nonzero elements of $k$.

	Suppose that $G$ is a finite group, $A$ is an abelian subgroup
of $G$, and that $(\rho, V)$ is an irreducible representation of $G$.
We can restrict $\rho$ to $A$ to get a representation of $A$ which may
or may not be irreducible.  At any rate there is a linear subspace $L$
of $V$ which is invariant under the restriction of $\rho$ to $A$, and
such that the restriction of $\rho$ to $A$ and to $L$ is an
irreducible representation of $A$.  It follows that $L$ is
$1$-dimensional, as in the preceding paragraph.  

	If $v$ is a nonzero vector in $L$, then $V$ is spanned by the
images of $v$ under $\rho_x$, $x \in G$, since $(\rho, V)$ is an
irreducible representation of $G$.  Let $E$ be a subset of $G$ such
that every element $x$ of $G$ can be expressed as $y \, a$ for some $y
\in E$ and $a \in A$, and so that the number of elements of $E$ is
equal to the number of elements of $G$ divided by the number of
elements of $A$.  In other words, $E$ should contain selections from
each of the cosets of $A$ in $G$.  Because $v$ is an eigenvalue of
$\rho_a$ for all $a \in A$, we obtain that the span of $\rho_x(v)$, $x
\in G$, which is the same as the span of $\rho_y(\rho_a(v))$ for $y
\in E$ and $a \in A$, is actually the same as the span of $\rho_y(v)$,
$y \in E$.  Therefore the dimension of $V$ is less than or equal to
the number of elements of $G$ divided by the number of elements of
$A$.

	Now let us assume that
\begin{equation}
	\hbox{$k$ is the field ${\bf C}$ of complex numbers}.
\end{equation}
Recall that a complex number $z$ can be written as $x + y \, i$, where
$x$, $y$ are real numbers, called the real and imaginary parts of $z$,
respectively.  The Fundamental Theorem of Algebra states that the
field of complex numbers is algebraically closed.

	If $z = x + y \, i$ is a complex number, with $x$, $y$ the
real and imaginary parts of $z$, then the complex conjugate of $z$ is
denoted $\overline{z}$ and defined to be $x - y \, i$.  If $z$, $w$
are complex numbers, then the complex conjugate of $z + w$ is the sum
of the complex conjugates of $z$ and $w$, and the complex conjugate of
$z \, w$ is the product of the complex conjugates of $z$ and $w$.  The
modulus of $z$ is denoted $|z|$ and is the nonnegative real number
such that $|z|^2 = x^2 + y^2$, which is the same as $z \,
\overline{z}$.  The triangle inequality states that $|z + w| \le |z| +
|w|$ for all complex numbers $z$, $w$.  One can also check that $|z \,
w| = |z| \, |w|$.

	Let $V$ be a finite-dimensional vector space over the complex
numbers.  By a Hermitian inner product on $V$ we mean a function
$\langle v, w \rangle$ defined for $v, w \in V$ and with values in the
complex numbers which satisfies the following properties.  First, for
each $w \in V$, $\langle v, w \rangle$ is a linear function in $v$.
Second, $\langle w, v \rangle$ is equal to the complex conjugate of
$\langle v, w \rangle$ for all $v, w \in V$.  Third, $\langle v, v
\rangle$ is a nonnegative real number for all $v \in V$ which is
equal to $0$ if and only if $v = 0$.

	Let $G$ be a finite group, and let $(\rho, V)$ be a
representation of $G$ over the complex numbers.  A Hermitian inner
product $\langle v, w \rangle$ on $V$ is said to be invariant under
the representation if $\langle \rho_x(v), \rho_x(w) \rangle$ is equal
to $\langle v, w \rangle$ for all $x \in G$ and $v, w \in V$.
If $\langle v, w \rangle_1$ is any Hermitian inner product on $V$,
then we can obtain a Hermitian inner product on $V$ from this one
which is invariant under the representation simply by summing
$\langle \rho_x(v), \rho_x(w) \langle_1$ over all $x \in G$.
Thus every representation of $G$ admits an invariant Hermitian
inner product.

	Let $G$ be a finite group, and let $(\rho, V)$ be a
representation of $G$ over the complex numbers which is equipped with
an invariant Hermitian inner product $\langle v, w \rangle$ on $V$.
Suppose that $L$ is a linear subspace of $V$ which is invariant under
$\rho$.  The orthogonal complement of $L$ in $V$ is denoted $L^\perp$
and consists of the vectors $v \in V$ such that $\langle v, w \rangle
= 0$ for all $w \in L$, and it is also invariant under the
representation since $L$ and the inner product are invariant.  In this
way one can decompose $(\rho, V)$ into an orthogonal direct sum of
irreducible representations.

	Let $G$ be a finite group, and let $(\rho, V)$ be a
representation of $G$.  If $x \in G$, then $\rho_x$ can be
diagonalized as a linear transformation on $V$, which is to say that
there is a basis of $V$ consisting of eigenvectors for $\rho_x$.  This
works just as well for an algebraically closed field $k$ of
characteristic $0$, because the subgroup of $G$ generated by $x$ is
abelian and the restriction of $\rho$ to this abelian subgroup can be
decomposed into a direct sum of $1$-dimensional representations.  In
the complex case one can argue instead that $\rho_x$ is a unitary
transformation with respect to an invariant inner product and hence
admits an orthonormal basis of eigenvectors.  Indeed, a unitary
transformation is normal, which is to say that it commutes with its
adjoint, and the existence of an orthonormal basis of eigenvectors can
be derived from the corresponding result for self-adjoint linear
transformations.

	Because $G$ is a finite group, $x^l$ is equal to the identity
element of the group for some positive integer $l$.  This implies that
$(\rho_x)^l$ is equal to the identity mapping on $V$, and therefore
the eigenvalues of $\rho_x$ are $l$th roots of unity.  The trace of
$\rho_x$, which is the same as the character of the representation
evaluated at $x \in G$, is therefore a sum of $l$th roots of unity,
and an algebraic integer in particular.

	A complex number which is an eigenvalue of $\rho_x$ has
modulus equal to $1$.  One can derive this from the fact that the
eigenvalues of $\rho_x$ are roots of unity, or using the fact that
$\rho_x$ is unitary with respect to an invariant inner product.  Hence
the inverse of an eigenvalue of $\rho_x$ is the same as its complex
conjugate.  The eigenvalues for $(\rho_x)^{-1}$ are the same as the
inverse of the eigenvalues for $\rho_x$, which are the complex
conjugates of the eigenvalues of $\rho_x$.  This can also be
seen in terms of $\rho_x$ being unitary, so that its inverse is
equal to its adjoint with respect to an invariant inner product.

	The trace of $(\rho_x)^{-1}$ is equal to the complex conjugate
of the trace of $\rho_x$.  This follows from the fact that $\rho_x$
and its inverse can be diagonalized, so that the trace is given by the
sum of the eigenvalues with their multiplicities.  One can also use
any orthonormal basis for the invariant inner product, since the
matrix for $(\rho_x)^{-1}$ using such a basis will be the adjoint of
the matrix for $\rho_x$, which is to say the complex conjugate of the
transpose of the matrix for $\rho_x$.  The diagonal entries for the
matrix for $(\rho_x)^{-1}$ are simply the complex conjugates of the
diagonal entries of the matrix for $\rho_x$.  If $\lambda$ is the
character associated to the representation, then it follows that
$\lambda(x^{-1})$ is equal to the complex conjugate of $\lambda(x)$
for all $x \in G$.

	As discussed before, $\lambda(x^{-1})$ is the same as the
character of the dual representation.  In other words, for a
representation of $G$ over the complex numbers, the character of the
dual representation is equal to the complex conjugate of the character
of the original representation.  One can also see this in terms of
matrices, if one describes a representation of $G$ on a vector space
of dimension $\ell$ in terms of a homomorphism from $G$ into the group
of invertible $\ell \times \ell$ matrices with complex entries.  The
existence of an invariant inner product amounts to being able to
describe the representation in terms of a homomorphism of $G$ into the
group of $\ell \times \ell$ unitary matrices, which are the matrices
with complex entries whose inverses are given by their adjoints or
conjugate transposes.  To get the dual representation one should take
the inverse transpose of the matrices, which is the same as the
complex conjugates of the matrices when they are unitary.

	In general, a representation of a group $G$ on a vector space
$V$ is a homomorphism from $G$ into the invertible linear
transformations on $V$, perhaps with additional regularity conditions,
such as continuity conditions.  One might wish to ask for some
additional data, like an invariant inner product on $V$ for a unitary
representation.  Many of the same notions as for representations of
finite groups on finite-dimensional vector spaces are applicable more
general, with elaborations as might be necessary.

\end{document}